%anomalous.tex: 
%%a Plain TeX file by Shalosh B. Ekhad (x pages)

%begin macros

\baselineskip=14pt
\parskip=10pt

\magnification=\magstephalf

\def\1{{\overline{1}}}
\def\2{{\overline{2}}}
\parindent=0pt
\overfullrule=0in

\def\frac#1#2{{#1 \over #2}}
%\headline={\rm  \ifodd\pageno  \RightHead  \else  \LeftHead  \fi}
%\def\RightHead{\centerline{
%Title
%}}
%\def\LeftHead{ \centerline{Doron Zeilberger}}
%end macros
\centerline
{\bf Automated Generation of Anomalous Cancellations}
\bigskip
\centerline
{\it Shalosh B. EKHAD} (with a Foreword by Doron Zeilberger)

{\bf Foreword} (By Doron Zeilberger)

In my last visit to Israel I bought a wonderful ``popular'' book [Sh], by the great expositor Haim Shapira, on
game theory. There he mentions (p. 156) how clueless students are. Once,  during a course on probability given in the Faculty of
Engineering, he solved a problem whose  answer happened to be
$$
\frac{16}{64} \quad .
$$
He `simplified' it, ``cancelling out" the the digit $6$:
$$
\frac{1\not6}{\not64} \, = \, \frac{1}{4} \quad ,
$$
getting the right answer. Shapira went on to say that he hopes that none of his readers reduces fractions this way,
and raised two questions, one to psychologists and one to mathematicians. The question to
psychologists was: {\it ``How come no (engineering!) student commented on this `simplification' for more than a minute a half?''}.
The question to mathematicians was: {\it ``Can you make up more  examples of this phenomenon?''.}

Not being a psychologist, I tried to focus on the second question.
When I got back home, I asked Neil Sloane whether he has seen this before, and, sure enough, he knew about this
phenomenon, and referred me to  sequence {\bf  A159975} and a few other related sequences in
the On Line Encyclopedia of Integer Sequences [Sl], that, in turn, referred me to a wonderful article [B]
by the great analyst (who  was also a great expositor!) Ralph Philip Boas.

In that article, Boas discusses the problem of finding such examples, in {\it any} given base.
He derived some complicated {\it non-linear} diophantine equations.

But why think so hard? If, given  a base $b$, and positive integers $d_1$ and $d_2$,
as well as the locations of the `cancelled' digits, $i_1$ and $i_2$, ($0 \leq i_1 <d_1$ and  $0 \leq i_2 <d_2$),
you are interested, in finding the set of {\it all}  fractions
$$
\frac{m}{n} \quad,
$$
written in base $b$, where $m$ and $n$ are with $d_1$ and $d_2$ digits, respectively, and where
the {\it illegal} cancellation of the $i_1$-th digit of the numerator and $i_2$-th digit of the denominator
gives you the same ratio, it suffices to generate all $(b^{d_1}- b^{d_1-1})\,(b^{d_2}- b^{d_2-1})$ such fractions and see
which ones work out. Alas, this takes {\it too long}, even nowadays, if you want to go anywhere far.

The compromise between being too clever and being too dumb, is to {\bf fix the denominator}, $n$,
and places $0 \leq i_1,i_2 \leq d_2-1$, and ask for a ``clever" way to find all numerators $m$ (of any size!) such that
if you `cancel out' the $i_1$-th digit of the top with the the $i_2$-th digit of the bottom, you get
the same thing. This leads to a {\it linear diophantine equation}, that is very easy to solve, thanks
to the {\it Extended Euclidean Algorithm}. Don Knuth ([K], p. 318, lines 5-6) calls the (unextended) Euclidean
algorithm the {\it granddaddy of all algorithms}, and then goes on to comment, that, most likely, it is not
due to Euclid (he only wrote it up in his textbook), and his `proof' of validity was not quite
up to his usual standards, since he lacked the  concept of {\it induction}.

The {\it Extended Euclidean Algorithm} is a simple by-product of the Euclidean algorithm, and in addition
to supplying the greatest common divisor of $A$ and $B$
$$
C \, = \, gcd(A,B) \quad,
$$
it also outputs integers $x$ and $y$ such that
$$
A\, x \, + \, B \, y  \, = \, C \quad .
$$

See [K], pp. 325-327, for a very nice and detailed account, in Knuth's  inimitable style.
According to an internet search, the Extended Euclidean Algorithm was first published by Roger Cotes (1682-1716)
who used it to compute  {\it continued fractions}.

So let's describe an algorithm that

{\bf inputs}

$\bullet$ A base $b$ \quad ;

$\bullet$ A positive integer $n$, written in base $b$, of $d_2$ digits, say (so $ b^{d_2-1} \leq n < b^{d_2}$) \quad ;

$\bullet$ Integers $i_1, i_2$ satisfying $0 \leq i_1,i_2 \leq d_2-1$ \quad ;

and

{\bf outputs}

$\bullet$ The set of all positive integers $m$ (of {\it any} size) such that if $m$ is written in base $b$,
the $i_1$-th digit of $m$ is the same as the $i_2$-th digit of $n$, and removing that common digit
from both $m$ and $n$ does not alter the ratio.

Since the denominator $n$ is assumed known, and we are singling out its $i_2$-th digit, let's call it $c$
($0 \leq c \leq b-1$), we can express it as
$$
n \, = \, N_1 \cdot b^{i_2+1} \, + \,  c \cdot b^{i_2}  \, +  \, N_2 \quad ,
\quad 1 \leq N_1 \leq b^{d_2-i_2-1} -1 \quad, \quad  0 \leq N_2 \leq b^{i_2} -1 \quad .
$$
We are  looking for integers $m$
$$
m \, = \, M_1 \cdot b^{i_1+1} \, + \,  c \cdot b^{i_1}  \, +  \, M_2 \quad ,
\quad 0 \leq M_1 \quad, \quad  0 \leq M_2 \leq b^{i_1} -1 \quad ,
$$
such that
$$
\frac{ M_1 \cdot b^{i_1+1} \, + \,  c \cdot b^{i_1}  \, +  \, M_2} { N_1 \cdot b^{i_2+1} \, + \,  c \cdot b^{i_2}  \, +  \, N_2} \, = \,
\frac{ M_1 \cdot b^{i_1} \,  +  \, M_2} { N_1 \cdot b^{i_2} \, +  \, N_2} \quad .
$$

Cross-multiplying, this leads to a certain {\it linear diophantine equation} of the form
$$
A \, M_1 \, + \, B  \, M_2 = C \quad,
$$
for some {\bf specific} integers $A,B,C$ (derivable from above), and unknowns $M_1,M_2$. If $C$ is not divisible by $gcd(A,B)$,
then there is no solution of course.
On the other hand, if $C$  is divisible by $gcd(A,B)$,  let $A_1 = \frac{A}{gcd(A,B)}$, $B_1 = \frac{B}{gcd(A,B)}$, $C_1 =\frac{C}{gcd(A,B)}$, and we get
$$
A_1 \, M_1 \, +  \, B_1 \, M_2 = C_1 \quad,
$$
where $gcd(A_1,B_1)=1$.
Using the Extended Euclidean algorithm, we manufacture
a pair of integers ($M_1^{(0)} ,M_2^{(0)}$) such that
$$
A_1 \, M_1^{(0)} \, + \, B_1 \, M_2^{(0)} = 1 \quad,
$$
from which we get the {\it general solution}
$$
M_1 \, = \,  C_1 \, M_1^{(0)} \, + \,  B_1 \, t\quad, \quad
M_2 \, = \, C_1 \, M_2^{(0)} \, - \, A_1 \, t .
$$
We then look for all $t$ that make $M_2$ between $0$ and $b^{i1}-1$.

{\bf The Maple package}

All this is implemented in the Maple package {\tt AnomalousCancellation.txt}, available directly from

{\tt http://sites.math.rutgers.edu/\~{}zeilberg/tokhniot/AnomalousCancellation.txt} \quad ,

or via the webpage of this article

{\tt http://sites.math.rutgers.edu/\~{}zeilberg/mamarim/mamarimhtml/anomalous.html} \quad ,

where there are several input and outputs files. For example, for a list of {\bf all} such anomalous cancellations
where the denominator is $\leq 9999$ (in our own base ten) see

{\tt http://sites.math.rutgers.edu/\~{}zeilberg/mamarim/mamarimPDF/ac4.pdf}  \quad ,

while for an extension of Boas' tables in pp. 119-122 of [B] see

{\tt http://sites.math.rutgers.edu/\~{}zeilberg/tokhniot/oAnomalousCancellation2.txt} \quad .

The above web-page has other output files, but  readers who have Maple can use the Maple package
{\tt AnomalousCancellation.txt} to get as many more examples as they wish.

\vfill\eject

{\bf References}

[B]  R. P. Boas, {\it  Anomalous Cancellation}, in:  Ch. 6 (pp. 113-129) in
{\it ``Mathematical Plums"}, Ross  Honsberger, ed., Dolciani Mathematical Expositions, Mathematical
Association of America, 1979.

[K] Donald E. Knuth, {\it ``The Art of Computer Programming'' vol. 2, Seminumerical algorithms}, 2nd ed.,
Addison-Wesley, 1981.

[Sh] Haim Shapira, ``{\it Conversations on Game Theory}'' (in Hebrew), Kinneret, Zmora-Bitan, Dvir, 2008.

[Sl] N. J. A. Sloane, {\it The On Line Encyclopedia of Integer Sequences}, {\tt  https://oeis.org} \quad .

\bigskip
\hrule
\bigskip
Shalosh B. Ekhad, c/o D. Zeilberger, Department of Mathematics, Rutgers University (New Brunswick), Hill Center-Busch Campus, 110 Frelinghuysen
Rd., Piscataway, NJ 08854-8019, USA. \hfill\break {\tt ShaloshBEkhad at gmail dot com}   \quad .

{\bf Exclusively published in the Personal Journal of Shalosh B. Ekhad and Doron Zeilberger and arxiv.org} \quad .

\end